%% file: root.tex
\documentclass[letterpaper, 10 pt, conference]{ieeeconf}  

\IEEEoverridecommandlockouts                              

\usepackage[noadjust]{cite}
\usepackage{graphicx} 
\graphicspath{ {./figures/} }
\usepackage{amsmath} 
\usepackage{threeparttable}
\usepackage{optidef}
\usepackage{multirow}\usepackage{threeparttable}
\usepackage{amssymb}
\usepackage{algorithm}
\usepackage[noend]{algpseudocode}
\usepackage{xcolor}
\definecolor{tumblue}{RGB}{0,101,189}

\usepackage{stfloats}
\usepackage{setspace}
\usepackage{tabularx}
\usepackage{booktabs}
\usepackage[T1]{fontenc}
\usepackage[scaled=0.85]{beramono}
\usepackage{listings}
\usepackage{caption}
\usepackage{cleveref}
\usepackage{chemformula}
\usepackage{url}
\usepackage{physics}
\Crefname{lstlisting}{Listing}{Listings}
\lstdefinestyle{sparqlstyle}{
  language=OCL,
  numbers=left,
  basicstyle=\footnotesize,
  stepnumber=1,
  numbersep=10pt,
  tabsize=2,
  showspaces=false,
  breaklines=true
}
\usepackage[nonumberlist,nopostdot,acronym,toc]{glossaries}
\usepackage{soul}
\usepackage{cases}

\usepackage{balance} 

\title{\LARGE \bf
Energy-Aware Model Predictive Control for Batch Manufacturing System Scheduling Under Different Electricity Pricing Strategies
}

\author{Hongliang Li$^{1}$, Herschel C. Pangborn$^{2}$, and Ilya Kovalenko$^{3}$
\thanks{$^{1}$Hongliang Li is with the Department of Industrial and Manufacturing Engineering, The Pennsylvania State University, University Park, PA, USA, e-mail: hjl5377@psu.edu.}%
\thanks{$^{2}$Herschel C. Pangborn is with the Department of Mechanical Engineering, The Pennsylvania State University, University Park, PA, USA, e-mail: hcpangborn@psu.edu.}%
\thanks{$^{3}$Ilya Kovalenko is with the  Department of Mechanical Engineering and the Department of Industrial and Manufacturing Engineering, The Pennsylvania State University, University Park, PA, USA, e-mail: iqk5135@psu.edu.}%
}

\begin{document}

\maketitle
\thispagestyle{empty}
\pagestyle{empty}

\begin{abstract}
Manufacturing industries are among the highest energy-consuming sectors, facing increasing pressure to reduce energy costs. This paper presents an energy-aware Model Predictive Control (MPC) framework to dynamically schedule manufacturing processes in response to time-varying electricity prices without compromising production goals or violating production constraints. A network-based manufacturing system model is developed to capture complex material flows, batch processing, and capacities of buffers and machines. The scheduling problem is formulated as a Mixed-Integer Quadratic Program (MIQP) that balances energy costs, buffer levels, and production requirements. A case study evaluates the proposed MPC framework under four industrial electricity pricing schemes. Numerical results demonstrate that the approach reduces energy usage expenses while satisfying production goals and adhering to production constraints. 
The findings highlight the importance of considering the detailed electricity cost structure in manufacturing scheduling decisions and provide practical insights for manufacturers when selecting among different electricity pricing strategies.

\end{abstract}

\section{Introduction}
\label{sec:intro}

\input{1_intro}
\section{Problem Statement}
\label{sec:problem}

\input{2_problem}
\section{Manufacturing System Model}
\label{sec:model}

\input{3_system_models_2}

\section{Model Predictive Control Formulation}
\label{sec:controller}

\input{4_MPC_controller_1}
\section{Case Study}
\label{sec:case}

\input{5_case_studies}
\section{Conclusion}
\label{sec:conclusion}

\input{6_conclusion}

\balance
\bibliographystyle{IEEEtran}
%

\bibliography{root}

\end{document}

%% file: 1_intro.tex

Modern manufacturing systems face increasing pressure to optimize schedules while balancing economic and environmental objectives.
The complex manufacturing configuration with batch processes and networked machines and buffers creates unique scheduling challenges~\cite{oliveira1993algorithms}.
Energy consumption is another critical factor, accounting for a significant share of total operational costs and contributing substantially to environmental impact~\cite{zhou2022iecl}.
In the U.S., the manufacturing sector is responsible for approximately 77\% of industrial electricity use, consuming around 27\% of the nation's total electricity~\cite{EIA2025}.
Recent advancements in smart grids and demand-side management offer manufacturers opportunities to reduce energy costs by controlling their production schedules in response to time-varying electricity prices~\cite{fadlullah2011toward,li2023system}.
However, manufacturers must ensure that these production shifts do not compromise the timely fulfillment of production goals to remain competitive.

Manufacturing scheduling has been extensively studied in operational research and control engineering.
Existing research formulates the scheduling problem as optimization programs~\cite{ pinedo2016scheduling,zhou2019real,shao2021ant,9737280,li2023system,xu2021industry}.
Some approaches focus on objectives such as minimizing makespan, reducing tardiness, and optimizing resource utilization, typically assuming static operating conditions and fixed resource costs~\cite{fang2011new, pinedo2016scheduling}.
Energy-aware manufacturing has emerged as an important research direction in response to rising energy costs.
Scheduling methods with energy considerations have been studied for machine shop systems~\cite{mouzon2007operational}.
Previous research demonstrated the potential for significant energy cost reduction through production scheduling that responds to variable electricity prices~\cite{shrouf2014optimizing}.
Recent work has expanded these approaches to multi-machine manufacturing environments, considering both energy and production efficiency~\cite{zhou2019real,shao2021ant}.
Existing methods are typically executed offline, limiting their ability to incorporate dynamic information and adapt to changing operational conditions.
Advancements in automation and the increasing integration of information technology in manufacturing present new opportunities for manufacturers to leverage highly dynamic information to enhance system responsiveness and efficiency~\cite{xu2021industry,li2023system}.
Therefore, there is a need for a scheduling framework that considers the complex manufacturing configurations as well as the dynamic nature of both production requirements and energy costs.

Model Predictive Control (MPC) has emerged as a powerful tool for systematically managing competing objectives and adapting to dynamic conditions while respecting constraints~\cite{qin1997overview}.
MPC operates by solving a finite-horizon optimization problem at each decision step using a predictive model to determine optimal control actions while continuously updating decisions based on new information~\cite{mayne2000constrained}.
MPC has been applied for manufacturing scheduling~\cite{baldea2014integrated,rubaiee2018energy}.
For single-machine scheduling, prior work incorporated energy considerations into an MPC framework under a non-preemptive scheduling scenario~\cite{rubaiee2018energy}.
To address more complex manufacturing environments with flexible processing sequences, MPC has also been integrated into multi-agent systems~\cite{9354443}.
Additionally, adapting manufacturing schedules to time-varying electricity prices has been an active area of research.
For instance,~\cite{GAGGERO2023845} demonstrated significant energy cost savings through MPC-based load shifting under Time-of-Use (TOU) pricing structures.

One important aspect of designing MPC for manufacturing systems is deriving an appropriate prediction model that considers the batch processing and complex manufacturing configurations.
Another important aspect is tuning the MPC optimization program to balance the trade-off between production requirements and energy costs.
While some studies have investigated MPC under specific time-varying pricing schemes, such as TOU pricing~\cite{GAGGERO2023845,li2023system}, a comprehensive assessment that considers a broader range of pricing structures remains unexplored in the literature.
This gap is particularly significant as industrial consumers face diverse electricity pricing options with complex cost structures, including both energy usage charge and additional demand charge based on peak power demand~\cite{ComEd_website, PJM_website}.

To address these limitations, we propose an MPC framework based on a network model of manufacturing systems.
This model captures parallel and serial production flows and batch processing.
A Mixed-Integer Quadratic Programming (MIQP) formulation balances time-varying energy costs, buffer levels, and production requirements while considering production constraints.
We provide a comprehensive evaluation of scheduling performance under four distinct electricity pricing programs commonly offered to industrial customers: General pricing, critical peak pricing, TOU pricing, and real-time pricing. By incorporating the cost structure of each pricing program, the evaluation provides practical insights into how each can be accounted for by manufacturing systems.
Through this work, we demonstrate how MPC-based scheduling can dynamically respond to electricity price signals while meeting production requirements in complex batch manufacturing environments.

The remainder of this paper is organized as follows.
\Cref{sec:problem} presents the assumptions and problem statement.
\Cref{sec:model} provides the manufacturing system model.
\Cref{sec:controller} formulates the energy-aware MPC.
\Cref{sec:case} presents case study results and discussions, while \Cref{sec:conclusion} provides conclusions and directions for future work.

%% file: 2_problem.tex
We consider a manufacturing system that simultaneously produces multiple distinct products through a network of interconnected processes and buffers.
To capture the networked configuration of the manufacturing systems as well as the batch processing, we propose a network-based model with continuous processing flows as shown in Fig.~\ref{fig:network_topology}.
Nodes represent buffers and edges represent manufacturing processes.
The nodes are categorized as intermediate buffers (storing the unfinished parts) and final buffers (storing the completed products).
The edges are categorized as intermediate edges (manufacturing processing) and terminal edges (product delivery). 
Each product has a production sequence with different processing steps.
The manufacturing system has batch processing that processes multiple parts simultaneously.
As shown in Fig.~\ref{fig:network_topology}, some of the processes are performed on the same machine, which need to be started and completed at the same time step, creating coupling constraints.

We make the following assumptions about the manufacturing systems:
(1) Each product has a daily production goal. These targets are specified by the higher-level planning team and assumed to be known a day ahead. 
(2) The system has unlimited raw materials for production.
(3) Electricity is the main energy source for the systems. Machine power is linearly related to the processing rate.
(4) The process yield rate is not considered.

Given this manufacturing environment, the central problem is to determine daily optimal production schedules considering time-varying electricity pricing that:
(1) meet the specified daily production goals,
(2) minimize the total energy cost incurred throughout the production process,
(3) maintain balanced material flows to avoid bottlenecks, buffer overflows, or production disruptions, and
(4) adhere to all production constraints.

%% file: 3_system_models_2.tex
\begin{figure}[t]
\centering
\includegraphics[width=1\linewidth]{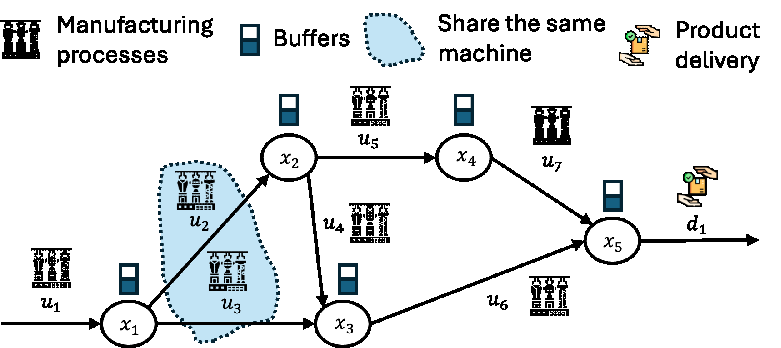}
\caption{Network topology of an example manufacturing system showing buffers (nodes), manufacturing processes (intermediate edges), and product delivery (terminal edge).}
\label{fig:network_topology}
\vspace{-15pt}
\end{figure}

\subsection{Manufacturing Network Model}
The discrete-time state-space representation of the system dynamics is:
\begin{equation}
\label{eq:sys-model}
x(k+1) = Ax(k) + Bu(k) + Ed(k)
\end{equation}
\noindent
At each time step $k$, $x(k) \in \mathbb{R}^{n_x}$ denotes the vector of buffer levels, $u(k) \in \mathbb{R}^{n_u}$ denotes the vector of processing flow rates, and $d(k) \in \mathbb{R}^{n_d}$ denotes the vector of final product outflows.
Each element $u_j(k)$ represents the processing rate along the edge, describing the rate at which material is transferred from its origin buffer to its destination buffer during time step $k$. Similarly, $x_i(k)$ denotes the inventory level in buffer $i$ at time $k$, and $d_p(k)$ captures the delivery quantity of final product $p$ at time $k$.

Matrix $A$ is the identity matrix, as the process yield rate is not considered. Matrix $B$ captures the network topology:
\begin{equation}
B_{ij} =
\begin{cases}
\hphantom{-}1, & \text{if $u_j$ flows into buffer $x_i$} \\
-1, & \text{if $u_j$ flows out of buffer $x_i$} \\
\hphantom{-}0, & \text{otherwise}
\end{cases}
\end{equation}
Matrix $E$ maps the final product outflows to their respective buffer states:
\begin{equation}
E_{ip} =
\begin{cases}
-1, & \text{if $d_p$ flows out of buffer $x_i$} \\
\hphantom{-}0, & \text{otherwise}
\end{cases}
\end{equation}
For numerical analysis, we define the outflow incidence matrix $B_o$ as:
\begin{align}
B_{o,ij} &=
\begin{cases}
-1, & \text{if } B_{ij} < 0 \\
\hphantom{-}0, & \text{otherwise}
\end{cases}
\end{align}
$B_o$ is defined to track material flows leaving each buffer in the manufacturing network.

The machine assignment matrix $M$ characterizes the mapping between manufacturing processing and available machines. Specifically:
\begin{equation}
M_{mj} =
\begin{cases}
1, & \text{if machine $m$ executes process $u_j$} \\
0, & \text{otherwise}
\end{cases}
\end{equation}
This binary mapping defines which machine is responsible for each process flow in the system.
Notably, certain machines may be assigned to multiple sequential or coupled processes, reflecting common manufacturing configurations. For example,  in Fig.~\ref{fig:network_topology}, processes $u_2$ and $u_3 $ are performed on the same machine.

To capture process coupling behavior, we define the process coupling matrix $C \in \{0,1\}^{n_u \times n_u}$ as:
\begin{equation}
C_{\varrho\rho} =
\begin{cases}
1, & \text{if $u_\varrho$ and $u_\rho$ are on the same machine} \\
0, & \text{otherwise}
\end{cases}
\end{equation}
The process coupling matrix $C$ identifies all pairs of processes that are physically coupled on a single machine. These processes must be activated and deactivated simultaneously at each time step.
\begin{figure}
    \centering
    \includegraphics[width=1\linewidth]{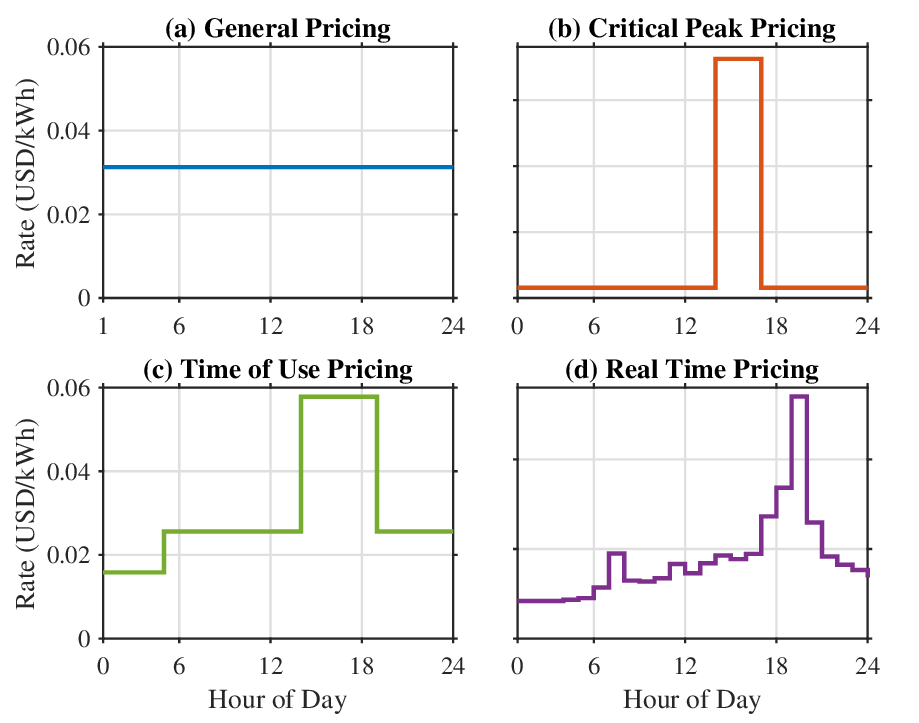}
    \caption{Chicago-area industrial electricity pricing programs in 2024: (a) General pricing, (b) Critical peak pricing, (c) TOU Pricing, and (d) Real-time pricing~\cite{ComEd_website,PJM_website}.}
    \label{fig:pricings}
    \vspace{-15pt}
\end{figure}
\subsection{Energy Consumption Model}
We characterize electricity pricing programs by considering three types of charges termed basic charge, demand charge, and energy usage charge~\cite{OpenEI_website,ComEd_website}.
The basic charge (USD/month) covers the fixed expenses of the utility company, such as transmission network investment and maintenance.
The demand charge (USD/kW) is based on the manufacturer's peak power load within a specified period, such as monthly.
The energy usage charge (USD/kWh) is based on electricity consumption.
For each unit of time, the energy usage of the system is calculated as:
\begin{equation}
\label{eq:energy_cal}
    E(k) = \sum_{j=1}^{n_u} \epsilon_j \cdot u_j(k)
\end{equation}
\noindent where
$u_j$ is the process rate (unit/hour), and $\epsilon_j$ (kWh/unit) is the power of the process.
Consider the discrete time step $k=1,\dots, N$, the total energy cost over the time horizon $N$ is calculated as:
\begin{equation}
    \begin{aligned}
        C_{total} &= C_{basic} + C_{demand} + C_{energy} \\
    & = \rho_b \cdot N + \rho_{d} \cdot max[E(k)] + \sum_{k =0}^{N-1}\rho_e(k)E(k)\\
    & \;\forall k=0,2,...,N-1
    \end{aligned}
\end{equation}
where
$C_{basic}$ is the basic charge calculated as the base charge rate $\rho_b$ multiplied by the horizon length,
$C_{demand}$ is the demand charge calculated as the demand charge rate $\rho_d$ multiplied by the peak power, and
$C_{energy}$ is the energy usage charge calculated as the time-varying energy cost $\rho_e(k)$ multiplied by the energy usage and summed over the horizon.
The demand charge calculation here is based on peak power within the time horizon of $N$, which may not be the true demand charges paid to the utility company.
In industry-standard practice, demand charges are calculated based on the peak power over 15 or 30-minute intervals throughout the billing period~\cite{ComEd_website}.



%% file: 4_MPC_controller_1.tex
This section presents the MPC formulation for runtime control of the manufacturing system, explicitly accounting for energy costs under different pricing policies.
The corresponding MIQP, along with its objective function and constraints, is detailed below.
\subsection{System Dynamics and State Input Constraints}
Equation~\eqref{eq:sys-model} serves as the dynamic prediction model.
The manufacturing system is subject to physical capacity limitations on all buffers and processes:
\begin{equation}
\label{eq: state input constraints}
    x_{\min} \leq x(k) \leq x_{\max},\;
    u_{\min} \leq u(k) \leq u_{\max}
\end{equation}
\noindent where
$x_{\min}$, $x_{\max}$, $u_{\min}$, and $u_{\max}$ are the minimum and maximum buffer and processing capacities, respectively.
Additionally, to ensure material conservation, outflows from any buffer must not exceed available material:
\begin{align}
\label{eq:outbuffer}
    Ax(k) + B_o u(k) \geq 0
\end{align}
Equation~\eqref{eq:outbuffer} ensures that buffer $x(k)$ has enough materials/parts for the next processing of $u(k)$.
Final product outflows are constrained by available material in the final product buffers:
\begin{align}
\label{eq:finalbuffer1}
    d(k) &\geq 0 \\
\label{eq:finalbuffer2}
    Ax(k) + Ed(k) &\geq 0
\end{align}
Equation~\eqref{eq:finalbuffer1} ensures that at each time step the production outflow is nonnegative, while
\eqref{eq:finalbuffer2} ensures that the final buffer has enough products for the delivery.
\subsection{Machine Activation Constraints}
Manufacturing schedules involve two key aspects: Determining when machines should be activated or deactivated and, once activated, specifying the processing rate.
We introduce binary variables $\delta_j(k) \in \{0,1\}$ to model the on/off status of each process to consider the discrete nature:
\begin{align}
\label{eq:off_condition}
    u_j(k) &\leq \Lambda \cdot \delta_j(k) \\
\label{eq:on_condition}
    u_j(k) &\geq \epsilon \cdot \delta_j(k)
\end{align}
where $\Lambda$ is a sufficiently large constant and $\epsilon$ is a small positive number representing the minimum flow rate when a process is active.
We identify all coupled process groups by the process coupling matrix \( C \).
Each coupled group \( g_i \subseteq \{u_1, u_2, \dots, u_{n_u}\} \) contains processes that are physically executed on the same machine.
For example, the system shown in Fig.~\ref{fig:network_topology} contains one coupled group: \( g_1 = \{u_2, u_3\} \).
To ensure that coupled processes are activated and deactivated simultaneously, we impose the following synchronization constraint on their binary activation variables:
\begin{equation}
\delta_{g_i[1]}(k) = \delta_{g_i[2]}(k) = \cdots = \delta_{g_i[|g_i|]}(k)
\end{equation}
where \( \delta_{g_i[j]}(k) \in \{0,1\} \) denotes the activation status of the \( j \)-th process in group \( g_i \) at time step \( k \).
This constraint ensures coordinated on/off decisions for all processes sharing the same machine.

\subsection{Production Requirement Constraints}
We formulate time-varying production requirement constraints that progressively tighten bounds as deadlines approach.
This formulation provides the MPC with more flexibility to shift production during earlier stages, while gradually prioritizing deadline compliance as the production horizon progresses.
Let $\pi_p(k)$ denote the predicted cumulative production of product $p$ at time step $k$. This is calculated as:
\begin{align}
\pi_p(k) = \sum_{t=0}^k d_p(k)
\end{align}
Let $\gamma_p$ denote the production goal of product $p$.
Let $\lambda_p$ denote the measured completed production at time $k=0$.
We introduce slack variables $s_p(k)$ to permit temporary deviations from production targets:
\begin{align}
\label{eq:slack_1}
\pi_p(k) &\geq \gamma_p-\lambda_p - s_p(k) \\
\label{eq:slack_2}
s_p(k) &\geq 0
\end{align}
Equation \eqref{eq:slack_1} ensures that the total production at time $k$ within the prediction horizon $N$ is near the remaining target by the deviation of $s_p(k)$. Equation
\eqref{eq:slack_2} ensures the nonnegative slack variables.
These soft constraints push the MPC controller to reach the production target diligently.

To avoid such deviation being far from the production target, we constrain the slack variable $s_p(k)$ with an upper bound defined by the production deviation allowance as:
\begin{align}
s_p(k) &\leq \alpha(k) \cdot \gamma_p
\end{align}
where $\alpha(k)$ is the allowable fraction of the production goals $\gamma_p$.
For example, if the production goal is $\gamma=100$ units with a production allowable fraction of $\alpha=0.05$, then the allowed deviation is $5$ units.
Instead of having the fixed $\alpha$, we progressively tighten $\alpha$ through the production.
We calculate the time-varying $\alpha(k)$ through a heuristic metric:
\begin{align}
\label{eq: production constraint 1}
\alpha(k) &= \tau \cdot \left(1 - \eta(k) \cdot (1 - \xi)\right) \\
\label{eq: production constraint 2}
\eta(k) &= \frac{h}{H} + \frac{k}{2N}
\end{align}
\noindent
where $\tau$ is the base production tolerance, $\xi$ is the tightening factor,
$\eta(k)$ is a progress metric that combines overall production progress with relative position within the prediction horizon,
$h$ is the time step of production, and $H$ is the length of the production horizon.
Note that $h$ and $H$ are distinct from $k$ and $N$, where the latter denote the time step of the MPC prediction and the length of the prediction horizon, respectively.



\subsection{Objective Function}
The MPC objective function balances buffer levels, production variation, slack penalty, and energy usage charges as:
\begin{equation}
\begin{aligned}
\label{eq:objective function}
\min_{u,d,\delta} \quad &J = \sum_{k=0}^{N-1} \Big[ \omega_x\|x(k)\|_2
+\omega_{u} \|u(k) - u(k-1)\|_2\\
&+ \omega_s \sum_{j=1}^{n_d}s_j(k) +\omega_e C_{energy} \Big]
\end{aligned}
\end{equation}
The first term penalizes the buffer level with weight $\omega_x$ to incentivize maintaining low buffer holding costs during production.
The second term penalizes production schedule variations with weight $\omega_u$ to ensure smooth production.
The third term penalizes the production slack with weight $\omega_s$.
The fourth term penalizes the energy usage charges during production with weight $\omega_e$.
Note that the demand and basic charges are not considered in the MPC formulation, as they relate to the long-term energy consumption profile.

\subsection{MPC Implementation Details}
The MPC controller dynamically updates machine schedules in a closed-loop manner. At each control step, the MPC solves the MIQP that minimizes the cost objective defined in~\eqref{eq:objective function}, subject to the manufacturing system dynamics in~\eqref{eq:sys-model} and constraints in~\eqref{eq: state input constraints}-\eqref{eq: production constraint 2}. 
A shrinking horizon implementation is used, where only the optimal decision variables at the first time step of each optimization are applied to the manufacturing system. This process repeats at every subsequent time step in a closed loop.

%% file: 5_case_studies.tex
\begin{figure}[t]
    \centering
    \includegraphics[width=0.9\linewidth]{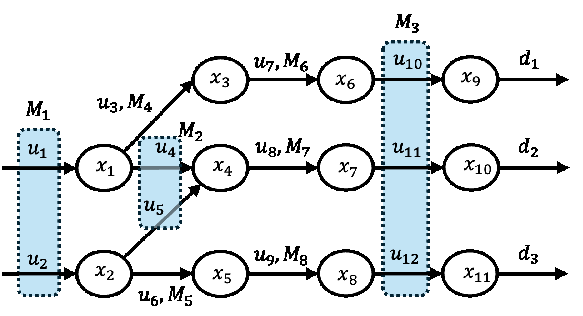}
    \caption{The manufacturing system considered in the case study includes 11 buffers and 12 processes, producing 3 products. Processes 1 and 2 share the same machine, as do processes 4, 5, and 10, 11, 12.}
    \label{fig:case_sys}
\end{figure}

\begin{table*}[t]
\renewcommand{\arraystretch}{1.3}
\caption{Industrial electricity pricing structures for Chicago area in 2024~\cite{ComEd_website, PJM_website}.}
\label{tab:pge_rates}
\setlength{\tabcolsep}{0.5em} 
\begin{tabular*}{\textwidth}{@{\extracolsep{\fill}}lccccc@{}}
\hline
\multirow{2}{*}{Electricity Pricing Program} 
    & \multirow{2}{*}{\shortstack{Basic Charge Rate\\(USD/Month)}} 
    & \multirow{2}{*}{\shortstack{Demand Charge Rate\\(USD/kW Peak Power/Month)}} 
    & \multicolumn{3}{c}{\shortstack{Energy Charge Rate (USD/kWh)}} \\
\cline{4-6}
 & & & Peak & Mid-Peak & Off-Peak \\
\hline
General Pricing
    & 221.77 
    & 10.93
    & ---
    & --- 
    & 0.03128 \\

Critical Peak Pricing$^1$ 
    & 221.77 
    & 10.93
    & 0.72500
    & --- 
    & 0.03128 \\

TOU Pricing$^2$ 
    & 221.77 
    & 10.93 
    & 0.05782 
    & 0.02561
    & 0.01583 \\

Real-Time Pricing 
    & 221.77 
    & 5.46 
    & \multicolumn{3}{c}{Hourly market rates} \\
\hline
\end{tabular*}

\vspace{1ex}
\raggedright 
$^1$ The critical peak period is from 1:00~pm to 5:00~pm on these 15 days. Avoiding critical peak events will generate 1 USD/day credit for the industrial customers.\\
$^2$ On-peak: 3:00~pm--8:00~pm weekdays; Mid-peak: 6:00~am--3:00~pm and 8:00~pm--10:00~pm weekdays; 
Off-peak: 10:00~pm--6:00~am all days, all hours weekends. \\
\vspace{-15pt}
\end{table*}

\begin{table}[t]
\renewcommand{\arraystretch}{1.3}
\caption{Simulation and MPC parameters.}
\label{tab:MPC_params}
\centering
\begin{tabular*}{\columnwidth}{@{\extracolsep{\fill}} l c l c }
\hline
Parameter & Value & Parameter & Value \\
\hline
$N$& 24 h & $\tau$, & 0.05 \\
$\xi$ & 0.5 & $\omega_s$ & 1000 \\
$\omega_x$ & 1 & $ \omega_u$ & 0.5 \\
$\omega_e$ & 500 & $ [\gamma_1,\gamma_2,\gamma_3]$ & [75,30,45] units\\
\hline
\end{tabular*}
\vspace{-15pt}
\end{table}

\begin{table}[t]
\renewcommand{\arraystretch}{1.3}
\caption{Manufacturing system parameters.}
\label{tab:Mfg_params}
\centering
\begin{tabular*}{\columnwidth}{@{\extracolsep{\fill}} l c l c }
\hline
Parameter & Value & Parameter & Value \\
\hline
$x^{1,\dots,11}_{min}$& 0 units & $x^{1,\dots,8}_{max}$, & 15 units \\
$x^{9,10,11}_{max}$ & 20 units & $u^{1,\dots,12}_{min}$& 0 units/hour\\
$u^{1,\dots,12}_{min}$& 0 units/hour& $u^{1,2}_{max}$ & 15 units/hour\\
$u^{3,\dots,12}_{max}$ & 10 units/hour& $\epsilon_{1,2}$ & 0.5 kWh/unit\\
$\epsilon_{10,11,12}$ & 0.5 kWh/unit & $\epsilon_{3,\dots,9}$ & 0.25 kWh/unit\\ 
\hline
\end{tabular*}
\vspace{-15pt}
\end{table}



\subsection{Case Study Setup}
We demonstrate the effectiveness of the proposed MPC framework with a realistic manufacturing system case study.
We evaluate the performance of the proposed method with one baseline general electricity pricing program and three time-varying electricity pricing programs, namely critical peak pricing, TOU pricing, and real-time pricing.
These four electricity pricing programs are described in Table~\ref{tab:pge_rates}. These pricing structures are representative of those
offered by major U.S. utility providers.
We consider a manufacturing system with 11 buffers, 12 processes, 8 machines, and 3 final products.
Fig.~\ref{fig:case_sys} shows the system network topology, where
$g_1 = \{u_1, u_2\}$ is assigned to machine 1, $g_2 = \{u_4, u_5\}$ is assigned to machine 2, and  $g_3 = \{u_{10}, u_{11}, u_{12}\}$ is assigned to machine 3.
The manufacturing scheduling problem is solved using the proposed MPC framework described in~\Cref{sec:controller}.
This is formulated in MATLAB~\cite{MATLAB} using YALMIP~\cite{Lofberg2004} and solved with the Gurobi optimization solver~\cite{gurobi}. The initial condition is set to 0 for all buffers.
For each pricing program, we simulate the manufacturing system over a 24-hour horizon with a known daily production goal of 75, 30, and 45 units for each product. The simulation parameters are summarized in Tables~\ref{tab:MPC_params}-\ref{tab:Mfg_params}.

\subsection{Results and Analysis}
\begin{figure}[t]
    \centering
    \includegraphics[width=1\linewidth]{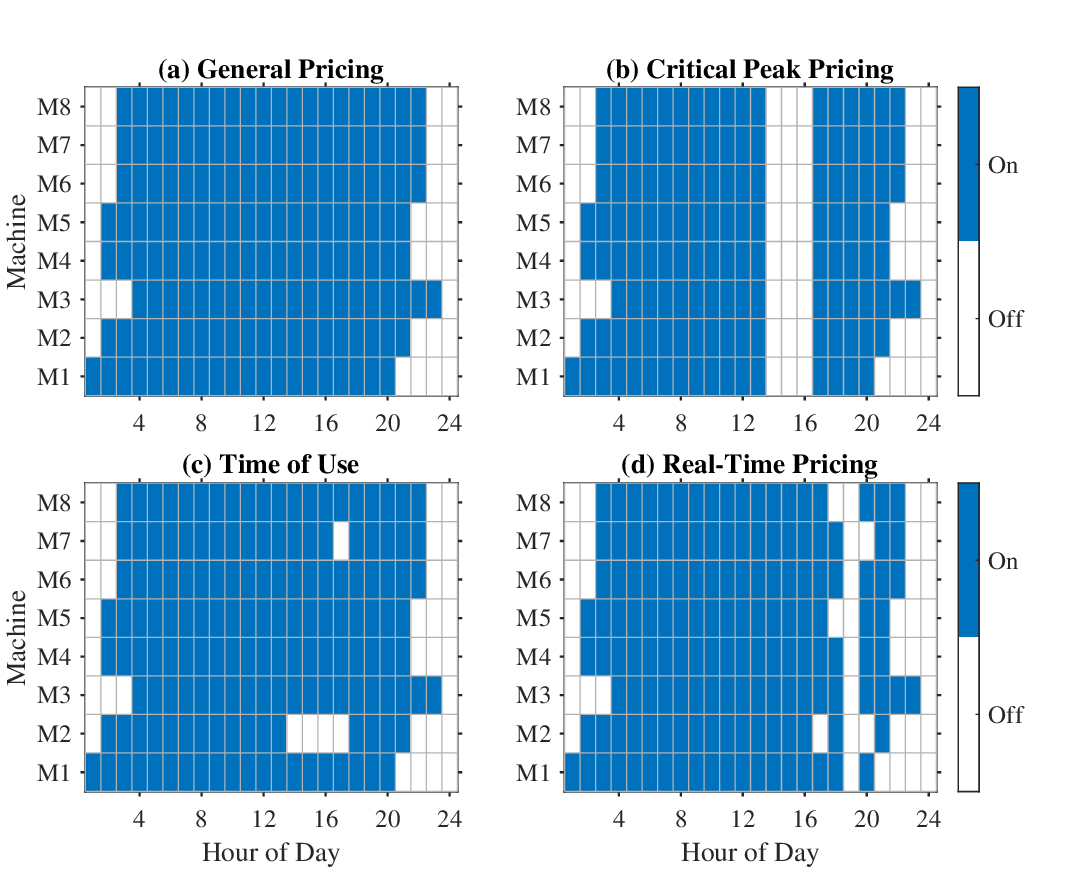}
    \caption{Machine on/off signal under different electricity pricing programs.}
    \label{fig:machine_on_off}
    \vspace{-14pt}
\end{figure}

Fig.~\ref{fig:machine_on_off} shows the machine activity patterns across all pricing programs over the 24-hour production horizon.
Machine activation patterns are effectively adapted to the underlying pricing structures.
Under general pricing, as shown in Fig.~\ref{fig:machine_on_off}(a), machines operate more uniformly throughout the day.
In contrast, dynamic pricing programs in Fig.~\ref{fig:machine_on_off}(b)-(d) show more strategic scheduling, concentrating operations during lower-price periods.
Under critical peak pricing, shown in Fig.~\ref{fig:machine_on_off}(b), there is an evident reduction in machine activity during the peak pricing event (hours 14-17). This demonstrates the controller's ability to respond to time-varying pricing by shifting production to lower-cost periods.

\begin{figure}[t]
    \centering
    \includegraphics[width=1\linewidth]{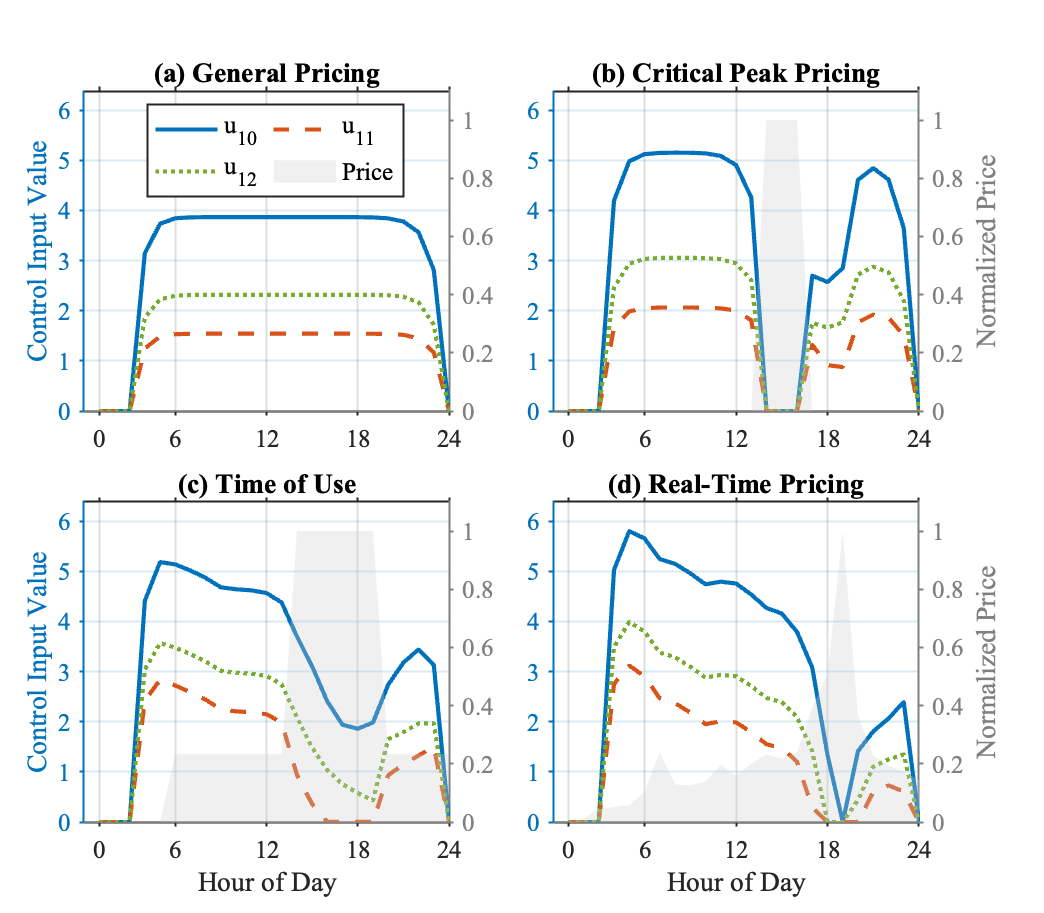}
    \caption{Processing rate of $u_{10}$, $ u_{11}$, and $ u_{12}$ under different electricity pricing programs.}
    \label{fig:process_rates}
    \vspace{-15pt}
\end{figure}

\begin{figure}[t]
    \centering
    \includegraphics[width=1\linewidth]{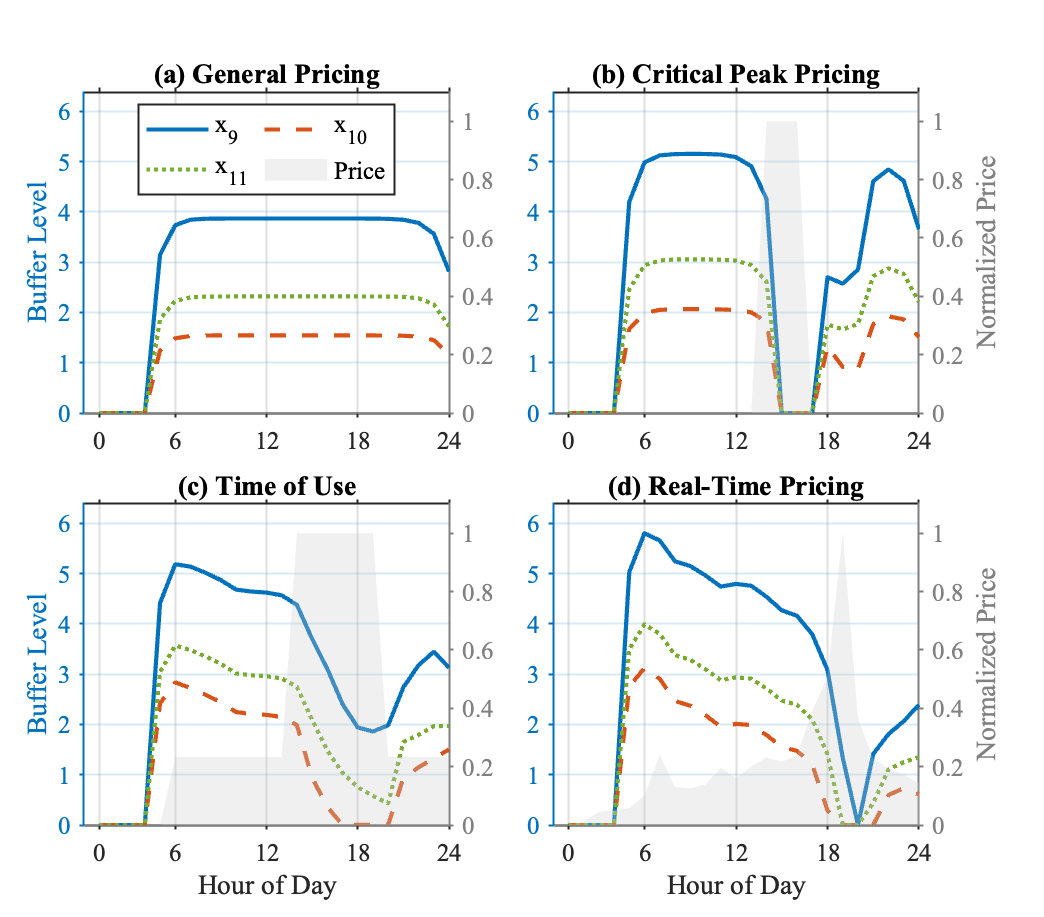}
    \caption{Buffer level of $x_9$, $ x_{10}$, and $ x_{11}$ under differenrt electricity pricing programs.}
    \label{fig:buffer_levels}
\end{figure}

\begin{figure}[t]
    \centering
    \includegraphics[width=1\linewidth]{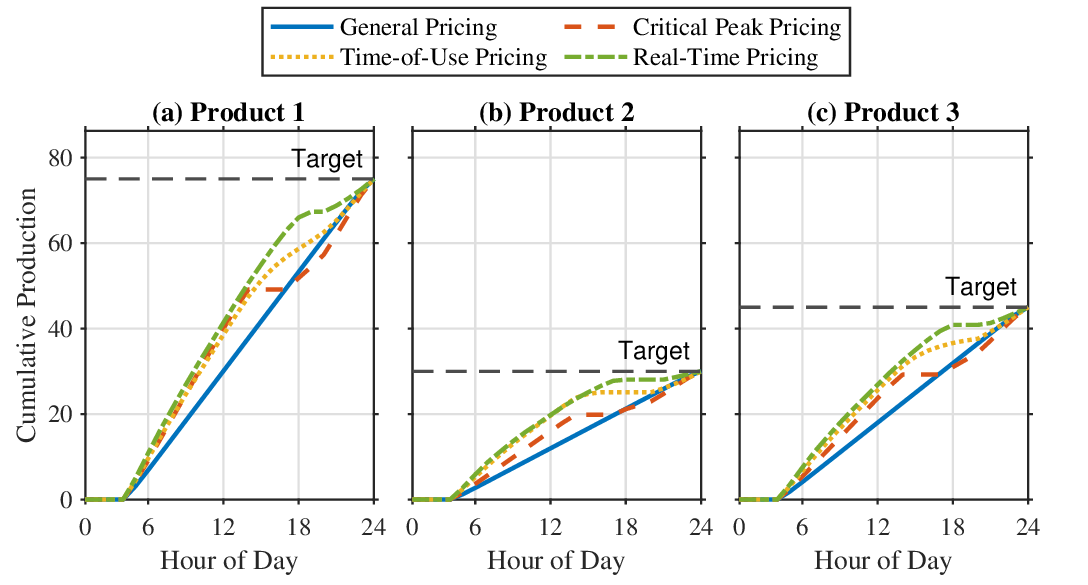}
    \caption{Cumulative production for each product under different pricing programs.}
    \label{fig:cumulative_production}
    \vspace{-15pt}
\end{figure}

Fig.~\ref{fig:process_rates} and Fig.~\ref{fig:buffer_levels} show how the MPC controller adjusts processing rates and manages inventory in response to different electricity pricing programs. 
Fig.~\ref{fig:process_rates} shows the processing rate of final processes $u_{10}$, $u_{11}$, and $ u_{12}$ across 24 hours.
Under general pricing with the uniform energy charge, the MPC algorithm generates consistent operations throughout most of the day.
Fig.~\ref{fig:process_rates}(a) exhibits minimal variation across the 24 hours, as the constant pricing structure incentivizes production scheduling based purely on meeting production requirements rather than the time-dependent energy cost.
The critical peak scenario displays an evident processing avoidance during the critical peak event period of hours 1:00 pm to 5:00 pm.
The processing rates show a significant reduction during these hours when electricity costs surge to $0.725$ USD/kWh, compared to the normal rate of $0.031$ USD/kWh. The system compensates by intensifying operations in periods immediately before and after the critical event.
TOU pricing produces different scheduling behavior with clear temporal patterns. Manufacturing activities concentrate during off-peak hours, with notable reductions from 3:00 pm to 8:00 pm. The system effectively shifts operations to leverage lower electricity rates during off-peak periods.
The real-time pricing scenario exhibits more dynamic and responsive scheduling.
This adaptive behavior is particularly evident in inputs $u_{10}$, $u_{11}$, and $u_{12}$, which display varying intensities throughout the day corresponding to hourly price fluctuations.
Fig.~\ref{fig:cumulative_production} shows the cumulative production for all three products.
The MPC controller drives the manufacturing system to successfully meet the daily production targets across different electricity pricing structures.
General pricing follows a relatively linear production pattern, while dynamic pricing models exhibit stepped approaches, with production accelerating during low-price periods. 

Table~\ref{tab:cost_comparison} presents the cost breakdown across all pricing programs, detailing basic charges, demand charges, energy usage charges, and total utility costs alongside peak power.
Analysis of energy usage charges reveals significant cost advantages for time-varying electricity pricing models.
Compared to general pricing, the critical peak, TOU, and real-time pricing programs achieved energy usage cost reductions of 14\%, 11\%, and 31\%, respectively.
However, these savings mechanisms involve strategic load shifting that concentrates operations during off-peak periods, consequently creating higher peak power under all dynamic pricing scenarios.
Notably, the real-time pricing model incorporates a substantially lower demand charge rate (5.46 USD/kW versus 10.93 USD/kW for other programs), mitigating the financial impact of its higher peak demand.
The real-time pricing rate structure represents a grid-level incentive for manufacturers to adopt more dynamic and responsive energy management strategies.
Real-time pricing achieves the lowest total cost at 109.02 USD, representing a 23\% reduction compared to general pricing. 
The MPC approach demonstrates adaptability in scheduling operations according to time-varying pricing while maintaining production targets. The controller successfully balances multiple objectives, including buffer management and energy cost minimization, across all pricing scenarios.
Critical peak, TOU, and real-time pricing programs effectively encourage shifting operations from high-cost periods to low-cost periods.
However, the load shifting creates higher peak power, which leads to higher demand charges.
For practical implementation, distributed energy resources, particularly on-site batteries or supercapacitors, could complement the MPC framework by mitigating peak power issues~\cite{garcia2018optimal}. 
This case study demonstrates that the relationship between electricity pricing programs and optimal manufacturing scheduling is complex. Intuitive approaches of shifting operations to lower-rate periods may not always yield the lowest total cost when demand charges constitute a significant portion of the electricity bill.

\begin{table*}[t]
\renewcommand{\arraystretch}{1.3}
\caption{Comparison of electricity cost under different pricing programs (all costs in USD).}
\label{tab:cost_comparison}
\centering
\begin{tabular*}{\textwidth}{@{\extracolsep{\fill}}lccccc@{}}
\hline
Pricing Programs& 
Basic Charge &
Demand Charge &
Peak Power (kW) &
Usage Charge &
Total Cost \\
\hline
General Pricing        & 7.39 & 126.73 & 11.59 & 7.04 & 141.16 \\
Critical Peak Pricing  & 7.39 & 168.27 & 15.39 & 6.04 & 181.70 \\
TOU Pricing            & 7.39 & 181.72 & 16.63 & 6.27 & 195.39 \\
Real-Time Pricing      & 7.39 & 96.79  & 17.73 & 4.84 & 109.02 \\
\hline
\end{tabular*}
\vspace{-10pt}
\end{table*}

%% file: 6_conclusion.tex
This paper presents an energy-aware Model Predictive Control (MPC) framework for manufacturing systems that optimally balances production requirements and energy cost minimization under dynamic electricity pricing.
The proposed approach models the manufacturing system as a network and leverages Mixed-Integer Quadratic Programming (MIQP) to make optimal scheduling decisions in response to time-varying electricity prices.
Case studies demonstrate that the proposed framework effectively shifts manufacturing operations away from high-price periods, leading to substantial reductions in energy usage costs while maintaining production targets.
However, the concentration of machine operations during lower-cost periods increases peak power demand, which may result in higher demand charges.
Future work can extend the proposed framework in two key directions: First is incorporating sources of uncertainty, such as volatile electricity prices and dynamic production goals, and second is integrating more complex, non-linear power profiles to improve real-world applicability.